\documentclass[preprint,3p,sort&compress,final,times]{elsarticle}
\usepackage{amsmath}
\usepackage{amssymb}
\usepackage{mathtools}
\usepackage{placeins}
\usepackage[usenames,dvipsnames]{color}
\usepackage{xcolor}
\usepackage{cases}
\usepackage{array}

\usepackage[usenames,dvipsnames]{color}

\journal{}
\begin{document}
\begin{frontmatter}

\title{Data driven non-equilibrium moist phase exchanges for atmospheric convection within a discontinuous Galerkin model 
of the compressible Euler equations}
\author[BOM]{David Lee\corref{cor}}
\ead{david.lee@bom.gov.au}
\author[ANU]{Kieran Ricardo}
\author[BOM]{Junwei Lyu}
\address[BOM]{Bureau of Meteorology, Melbourne, Australia}
\address[ANU]{Mathematical Sciences Institute, Australian National University, Canberra, Australia}
\cortext[cor]{Corresponding author. Tel. +61 452 262 804.}

\begin{abstract}
A discrete formulation of the multiphase moist compressible Euler equations is presented for which the energy and tracer variance 
conserving dynamics is coupled to the non-equilibrium moist phase exchanges in a thermodynamically consistent manner.
A neural network is then trained to learn the mass exchanges between vapour, liquid and ice phases in atmospheric convection. The network is 
trained on convection resolving output from a regional configuration of the LFRic model with a multi-moment microphysics parameterisation (CASIM). 
The loss function for learning these phase exchanges is formulated under the assumptions of thermal and mechanical equilibrium (same temperature 
and pressure for all phases), and chemical dis-equilibrium (different Gibbs free energies for all phases). The network outputs determine the 
exchanges of vapour, liquid and ice so as to conserve mass, and the resulting change in entropy is determined from the network outputs so as to 
conserve energy. The neural network is implemented in a thermodynamically consistent manner within 
the discontinuous Galerkin moist compressible Euler model
in order to simulate the formation of three-phase clouds for convection at sub-km resolution.
The results are compared to those from a physics based representation of three-phase moist processes at thermodynamic equilibrium, as well as
to an idealised configuration of LFRic model.
\end{abstract}

\end{frontmatter}

Convection resolving atmospheric models require sub-grid parameterisations for micro-physical processes involving phase exchanges 
between different moisture species. Sophisticated multi-moment parameterisations of these processes can be computationally expensive, since they 
require the advection of both the mass and number fractions for each moisture species (and potentially other moments as well) 
\cite{Shipway12, Field23}. These parameterisations are typically implemented within atmospheric models under the assumptions of mechanical 
and thermal equilibrium (a single temperature and pressure for all phases), and chemical dis-equilibrium (slow variations in mass fractions 
between species resulting in different Gibbs free energies for all species) \cite{BowenThesis}.

In the present work we seek to learn the mass exchanges between vapour, liquid and ice for convection resolving, non-hydrostatic atmospheric 
dynamics at chemical dis-equilibrium via output from a regional configuration of the LFRic model \cite{Johnson26}, and to implement the resulting 
neural network within an idealised model of the moist compressible Euler equations in order to simulate phase changes in high resolution convective 
dynamics. The motivation for doing so is two-fold. Firstly, the inference of a neural network for determining phase changes may be more 
computationally efficient than the existing physics based scheme, since it negates the need to explicitly advect additional moments. Secondly, the 
abstraction of the phase exchanges via a data driven scheme allows us to formulate these exchanges in a thermodynamically consistent manner, with 
the intrinsic thermodynamic variables (pressure, temperature and chemical potentials for each phase) derived from the corresponding prognostic 
variables (density, entropy, mass fractions for each phase) via a thermodynamic potential (in this case the internal energy). Moreover the 
thermodynamically consistent formulation allows us to derive the resulting change in entropy (or other thermal variable) due to phase exchanges so 
as to exactly conserve energy.

Thermodynamically consistent formulations of moist atmospheric processes have been previously implemented in idealised numerical models for both 
equilibrium \cite{Thuburn17,Bowen22a,Ricardo26} and non-equilibrium \cite{Bowen22b} configurations. In these equilibrium formulations, a Newton problem
is typically solved at each time step (or sub-step) to determine the partition of the moisture (and sometimes entropy) between phases by matching
the Gibbs free energies (and temperatures) of the different phases. In the non-equilibrium case, a coupled system of prognostic equations is solved
for the mass fractions and entropies of the different phases, subject to algebraic constraints to ensure the positivity of the mass of the different 
phases. In the present data-driven approach, only chemical dis-equilibrium is assumed, and when applied in the dynamical model, the thermodynamics
is embedded in the neural network, such that a coupled minimisation problem is not required. Thermodynamically consistent machine learning techniques 
have been applied
previously to learn the manifold of thermodynamic potentials for other systems such as the Van der Waals equation of state and for a Lennard-Jones 
system \cite{Rosenberger22}, with the derivatives of the thermodynamic potential embedded in the loss function.
There have also been numerous implementations of neural network substitutes for moist physics parameterisations within global atmospheric models
\cite{Brenowitz18,Han20,Wang22,WattMeyer24}, with a focus more on directly replicating the existing physical parameterisations than on considerations 
of thermodynamic consistency.

The remainder of this article proceeds as follows: In Section 1 a representation of non-equilibrium thermodynamics within an energy conserving
dynamical formulation of the moist compressible Euler equations is presented. In Section 2 we describe the training of a neural network to represent
the mass exchanges between vapour, liquid and ice, and to derive the resulting entropy change subject to conservation of energy. The discontinuous
Galerkin numerical model of the three phase moist compressible Euler equations into which the neural network is integrated is detailed in Section 3.
Results for the simulation of a standard test case for a moist bubble and comparisons to an existing representation of moist processes at thermodynamic
equilibrium, and to the original LFRic model with CASIM microphysics are presented in Section 4. Finally some conclusions are discussed in Section 5.

\section{Conservative formulation of non-equilibrium moist processes for the compressible Euler equations}

Energy conserving Hamiltonian systems, including non-canonical systems such as the compressible Euler equations 
\cite{Lee21,Cotter23} may be presented as
\begin{equation}
\frac{\partial y}{\partial t} = \boldsymbol{\mathsf{A}}\frac{\delta\mathcal{H}}{\delta y},
\end{equation}
where $y$ constitutes a vector of prognostic variables, $\mathcal{H}(y)$ represents the energy and 
$\boldsymbol{\mathsf{A}}(y)=-\boldsymbol{\mathsf{A}}^{\top}(y)$ represents some skew-symmetric matrix 
\cite{Hairer06}. Energy is conserved by construction for such systems as 
\begin{equation}
\frac{\mathrm{d}\mathcal{H}}{\mathrm{d}t} = \Bigg(\frac{\delta\mathcal{H}}{\delta y}\Bigg)^{\top}\frac{\partial y}{\partial t} = 
\Bigg(\frac{\delta\mathcal{H}}{\delta y}\Bigg)^{\top}\boldsymbol{\mathsf{A}}\frac{\delta\mathcal{H}}{\delta y} = 0
\end{equation}
due to the skew-symmetry of $\boldsymbol{\mathsf{A}}$. The above formulation accounts for energy conserving dynamical processes. 
Irreversible thermodynamic source terms may be incorporated into such a formulation via the introduction of a symmetric operator
$\boldsymbol{\mathsf{B}}=\boldsymbol{\mathsf{B}}^{\top}$ as \cite{Eldred20}
\begin{equation}\label{hamiltonian_AB}
\frac{\partial y}{\partial t} = \boldsymbol{\mathsf{A}}\frac{\delta\mathcal{H}}{\delta y} + 
\boldsymbol{\mathsf{B}}\frac{\delta\mathcal{H}}{\delta y}.
\end{equation}

For the moist compressible Euler equations under the assumptions of mechanical and thermal equilibrium (constant pressure, $p$ and temperature, 
$T$ for all phases) and chemical dis-equilibrium (different Gibbs free energies for all phases) \cite{BowenThesis}, the energy over the domain 
$\Omega$ may be defined with respect to the
prognostic variables of wind $\boldsymbol{u}$, density $\rho$, entropy $\eta$ and mass fractions of vapour $q_v$, liquid $q_l$ and ice $q_i$ as
\begin{equation}\label{energy}
\mathcal{H} = \int\frac{1}{2}\rho\boldsymbol{u}^2 + \rho gz + \rho u(\rho,\eta,q_d,q_v,q_l,q_i)\mathrm{d}\Omega,
\end{equation}
where $g$ is the gravitational acceleration, $z$ is the vertical coordinate, $q_d = 1 - q_v - q_l - q_i$ is the mass fraction of dry air and 
$u(\rho,\eta,q_d,q_v,q_l,q_i)$ is the specific internal energy. In order to enforce thermodynamic consistency, the intrinsic thermodynamic variables
(pressure $p$, temperature $T$ and the chemical potentials of the different phases $\mu_j$, $j\in\{d,v,l,i\}$) may be determined as the variational
derivatives of the internal energy with respect to the prognostic variables as \cite{Eldred22}
\begin{subequations}\label{du}
\begin{align}
p &= \frac{1}{\rho^2}\frac{\delta u}{\delta\rho} \\
T &= \frac{\delta u}{\delta\eta} \\
\mu_j &= \frac{\delta u}{\delta q_j}.
\end{align}
\end{subequations}
Importantly for our energy conserving formulation, the variational derivative of the energy with respect to the density is given as:
\begin{equation}\label{dHdrho}
\frac{\delta\mathcal{H}}{\delta\rho} = \frac{1}{2}\boldsymbol{u}^2 + gz + u + \frac{p}{\rho} := \Phi.
\end{equation}

The momentum equation for the compressible Euler equations in vector invariant form subject to a gravitational body force is most commonly expressed as
\begin{equation}\label{momentum}
\frac{\partial\boldsymbol{u}}{\partial t} + (\nabla\times\boldsymbol{u})\times\boldsymbol{u} + \frac{1}{2}\boldsymbol{u}^2 + gz + \frac{1}{\rho}\nabla p = 0.
\end{equation}
Recalling the identity $\nabla u = \frac{\delta u}{\delta\rho}\nabla\rho + \frac{\delta u}{\delta\eta}\nabla\eta + \sum_j\frac{\delta u}{\delta q_j}\nabla q_j$
\cite{Bannon03} and \eqref{du}, \eqref{dHdrho}, we may re-express the pressure gradient term in \eqref{momentum} as
\begin{equation}
\frac{1}{\rho}\nabla p = \nabla\Big(u + \frac{p}{\rho}\Big) - T\nabla\eta - \sum_j\mu_j\nabla q_j.
\end{equation}
We may then cast the moist compressible Euler equations in chemical dis-equilibrium with respect to the energy \eqref{energy} in the form presented
in \eqref{hamiltonian_AB} as
\begin{equation}\label{euler_d}
\begin{bmatrix}
\frac{\partial \boldsymbol{u}}{\partial t}\\ \frac{\partial \rho}{\partial t} \\ \frac{\partial \eta}{\partial t} \\ 
\frac{\partial q_d}{\partial t} \\ \frac{\partial q_v}{\partial t} \\ \frac{\partial q_l}{\partial t} \\ \frac{\partial q_i}{\partial t}
\end{bmatrix} = 
\begin{bmatrix}
-\frac{\nabla\times\boldsymbol{u}}{\rho} & -\nabla & \frac{\nabla\eta}{\rho} & \frac{\nabla q_d}{\rho} & \frac{\nabla q_v}{\rho} & \frac{\nabla q_l}{\rho} & \frac{\nabla q_i}{\rho} \\
-\nabla\cdot & 0 & 0 & 0 & 0 & 0 & 0 \\
-\frac{\nabla\eta}{\rho} & 0 & 0 & 0 & 0 & 0 & 0 \\
-\frac{\nabla q_d}{\rho} & 0 & 0 & 0 & 0 & 0 & 0 \\
-\frac{\nabla q_v}{\rho} & 0 & 0 & 0 & 0 & 0 & 0 \\
-\frac{\nabla q_l}{\rho} & 0 & 0 & 0 & 0 & 0 & 0 \\
-\frac{\nabla q_i}{\rho} & 0 & 0 & 0 & 0 & 0 & 0 
\end{bmatrix}
\begin{bmatrix}
\rho\boldsymbol{u} \\ \Phi \\ \rho T \\ \rho\mu_d \\ \rho\mu_v \\ \rho\mu_l \\ \rho\mu_i
\end{bmatrix} + \\
\begin{bmatrix}
0 & 0 & 0 & 0 & 0 & 0 & 0 \\
0 & 0 & 0 & 0 & 0 & 0 & 0 \\
0 & 0 & A & 0 & 0 & 0 & 0 \\
0 & 0 & 0 & 0 & 0 & 0 & 0 \\
0 & 0 & 0 & 0 & {B}+{C} & -{B} & -{C} \\
0 & 0 & 0 & 0 & -{B} & {B}+{D} & -{D} \\
0 & 0 & 0 & 0 & -{C} & -{D} & {C}+{D} 
\end{bmatrix}
\begin{bmatrix}
\rho\boldsymbol{u} \\ \Phi \\ \rho T \\ \rho\mu_d \\ \rho\mu_v \\ \rho\mu_l \\ \rho\mu_i
\end{bmatrix},
\end{equation}
where $B$, $C$, $D$ represent the mass exchanges from vapour to liquid, vapour to ice, and liquid to ice respectively, and $A$ represents 
the change in entropy due to phase exchanges. 
The first matrix in \eqref{euler_d} is the skew-symmetric operator accounting for the dynamics (transport, pressure gradients, and vorticity), 
and conserves energy by construction since left multiplication of this term by the vector $\delta\mathcal{H}/\delta y$ for 
$y=[\boldsymbol{u},\rho,\eta,q_d,q_v,q_l,q_i]$ vanishes. The
second term meanwhile represents the symmetric operator accounting for the source terms due to thermodynamic processes.
In order to conserve energy for the thermodynamic source terms, the change in entropy due to phase exchanges may be 
determined by left multiplication of this term by $\delta\mathcal{H}/\delta y$ such that
\begin{equation}\label{entropy_inc}
	A = -\frac{1}{T^2}\Big((\mu_v-\mu_l)^2{B} + (\mu_v-\mu_i)^2{C} + (\mu_l-\mu_i)^2{D}\Big).
\end{equation}

Given that the dry mass fraction may be inferred directly from the sum of the moist mass 
fractions as $q_d = 1 - \sum_kq_k$, $k\in\{v,l,i\}$, we may eliminate the advection equation for $q_d$ and incorporate the dry air pressure 
gradient forcing term $\mu_d\nabla q_d = -\mu_d\sum_k\nabla q_k$ directly into those for the other mass fractions as $\mu_k' = \mu_k - \mu_d$, 
resulting in a reduced system for the moist compressible Euler equations as
\begin{equation}\label{euler}
\begin{bmatrix}
\frac{\partial \boldsymbol{u}}{\partial t}\\ \frac{\partial \rho}{\partial t} \\ \frac{\partial \eta}{\partial t} \\ 
\frac{\partial q_v}{\partial t} \\ \frac{\partial q_l}{\partial t} \\ \frac{\partial q_i}{\partial t}
\end{bmatrix} = 
\begin{bmatrix}
-\frac{\nabla\times\boldsymbol{u}}{\rho} & -\nabla & \frac{\nabla\eta}{\rho} & \frac{\nabla q_v}{\rho} & \frac{\nabla q_l}{\rho} & \frac{\nabla q_i}{\rho} \\
-\nabla\cdot & 0 & 0 & 0 & 0 & 0 \\
-\frac{\nabla\eta}{\rho} & 0 & 0 & 0 & 0 & 0 \\
-\frac{\nabla q_v}{\rho} & 0 & 0 & 0 & 0 & 0 \\
-\frac{\nabla q_l}{\rho} & 0 & 0 & 0 & 0 & 0 \\
-\frac{\nabla q_i}{\rho} & 0 & 0 & 0 & 0 & 0 
\end{bmatrix}
\begin{bmatrix}
\rho\boldsymbol{u} \\ \Phi \\ \rho T \\ \rho\mu_v' \\ \rho\mu_l' \\ \rho\mu_i'
\end{bmatrix} + \\
\begin{bmatrix}
0 & 0 & 0 & 0 & 0 & 0 \\
0 & 0 & 0 & 0 & 0 & 0 \\
0 & 0 & A & 0 & 0 & 0 \\
0 & 0 & 0 & {B}+{C} & -{B} & -{C} \\
0 & 0 & 0 & -{B} & {B}+{D} & -{D} \\
0 & 0 & 0 & -{C} & -{D} & {C}+{D} 
\end{bmatrix}
\begin{bmatrix}
\rho\boldsymbol{u} \\ \Phi \\ \rho T \\ \rho\mu_v \\ \rho\mu_l \\ \rho\mu_i
\end{bmatrix}.
\end{equation}
Energy is still conserved for the skew-symmetric form with the elimination of $q_d$, since the temporal chain rule still holds since
$\sum_k\mu_k'\frac{\partial q_k}{\partial t} = \sum_k\mu_k\frac{\partial q_k}{\partial t} - \mu_d\sum_k\frac{\partial q_k}{\partial t} = 
\sum_k\mu_k\frac{\partial q_k}{\partial t} + \mu_d\frac{\partial q_d}{\partial t}$.

\section{Training of the neural network}

The formulation of the moist compressible Euler equations given in \eqref{euler} gives no details as to the physical description of the
phase exchanges between vapour, liquid and ice encoded in the operators $B, C, D$. Here these are learned via a neural network of the 
form $\mathbb{R}^5\rightarrow\mathbb{R}^3: (\rho,\eta,q_v,q_l,q_i)\rightarrow(\tilde{B},\tilde{C},\tilde{D})$ which takes as its inputs 
the thermodynamic state determined from the prognostic variables, and returns the weighted phase exchanges. In order to learn these 
exchanges we define the loss function via the transport equations for the mass fractions as given in \eqref{euler} as
\begin{equation}\label{loss_func}
\mathcal{L} := \mathcal{L}_v + \mathcal{L}_l + \mathcal{L}_i \approx 0,
\end{equation}
where
\begin{subequations}\label{loss_func_2}
\begin{align}
\mathcal{L}_v &= \Bigg(\frac{\partial q_v}{\partial t} + \boldsymbol{u}\cdot\nabla q_v - 
\rho (q_v+q_l)\tilde{B}(\mu_v-\mu_l) - \rho (q_v+q_i)\tilde{C}(\mu_v-\mu_i)\Bigg)^2, \\
\mathcal{L}_l &= \Bigg(\frac{\partial q_l}{\partial t} + \boldsymbol{u}\cdot\nabla q_l - 
\rho (q_v+q_l)\tilde{B}(\mu_l-\mu_v) - \rho (q_l+q_i)\tilde{D}(\mu_l-\mu_i)\Bigg)^2, \\
\mathcal{L}_i &= \Bigg(\frac{\partial q_i}{\partial t} + \boldsymbol{u}\cdot\nabla q_i - 
\rho (q_v+q_i)\tilde{C}(\mu_i-\mu_v) - \rho (q_l+q_i)\tilde{D}(\mu_i-\mu_l)\Bigg)^2.
\end{align}
\end{subequations}
Since the chemical potentials $\mu_k$ are many orders of magnitude larger than the mass fractions, $q_k$, we 
scale the transfer terms as $B=(q_v+q_l)\tilde{B}$, $C=(q_v+q_i)\tilde{C}$, $D=(q_l+q_i)\tilde{D}$ in order to improve 
convergence of the network training.

The network is trained on data from a 6 hour run of a regional configuration of the LFRic model \cite{Johnson26},
centered on Darwin, Australia from 3.30pm to 9.30pm local time on the $6^{th}$ of January 2021 (a time period of 
strong tropical convection), at 1.5km resolution, with a 60s time step. The model uses the double-moment configuration of the
CASIM microphysics 
parameterisation \cite{Shipway12,Field23}, and is run on a domain of $480\times 480$ cells with 70 vertical levels
and a domain top of 30km. 
The model data is taken from the second 6-hourly cycle to ensure that the second moment of the CASIM parameterisation 
(the number concentration) was spun up to a physically realistic distribution. This is despite the fact that the number 
concentration is not used as an input for training the network. 
The model data was filtered by eliminating the top and bottom two model levels, and also the 120 
cells on each of the horizontal boundaries to eliminate potential artefacts from the imposition of lateral boundary 
conditions. The data was also filtered for samples where the value of the mass fraction $q_k$ is greater than
$10^{-7}$ for at least one water phase, and every second model time step was used owing to the large size of the output files.

The spatial and temporal derivatives of the material advection operator in the loss function \eqref{loss_func_2} were
reconstructed from the model output data using centered second order finite differences in space and time. The liquid mass 
fraction was determined as the sum of the cloud liquid and rain mass fractions, and the ice mass fraction
was given as the sum of the cloud ice, snow and graupel mass fractions, while the density and pressure were derived
from the native LFRic model variables of potential temperature and Exner pressure.
See Appendix A for further details on the data conversion.

The network is of the form 
$\mathbb{R}^5\rightarrow\mathbb{R}^{10}\rightarrow\mathbb{R}^{60}\rightarrow\mathbb{R}^{60}\rightarrow\mathbb{R}^{60}\rightarrow\mathbb{R}^{12}\rightarrow\mathbb{R}^3$, 
and uses leaky ReLU activation functions at each level, with an Adam optimiser. The initial learning rate 
is set as $1.4\times 10^{-3}$, and is subsequently reduced to $7\times 10^{-4}$ and then $3.5\times 10^{-4}$.
The loss is computed from the mean square error of the loss function in \eqref{loss_func}. At 
each epoch the data is batched to a size of 16,384 samples. 

\begin{figure}[!hbtp]
\centering
\includegraphics[width=0.60\textwidth,height=0.45\textwidth]{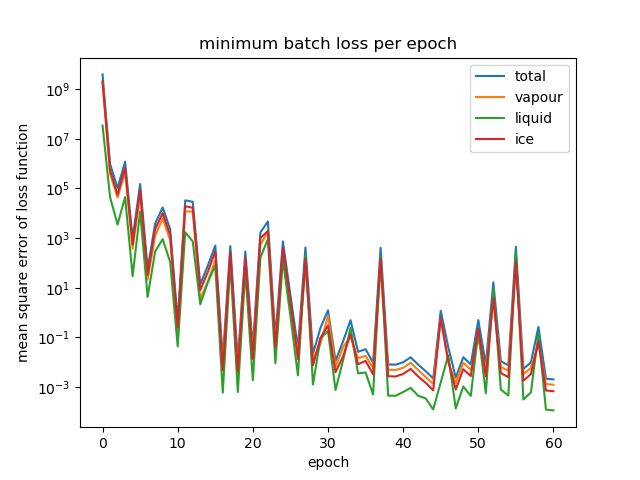}
\caption{Minimum batch loss per epochs for the total loss function, as well as for the individual
vapour, liquid and ice components.}
\end{figure}

Figure 1 shows the minimum batch loss per epoch for the total loss function as well as for the individual vapour, liquid and 
ice components. The final total loss is $2\times 10^{-12}$ the size of the initial loss.
Training to such a small loss magnitude is necessary owing to the large disparity in the
magnitudes of the left and right hand sides of the loss function.

Note that while it is possible to reduce the loss futher, experiments of moist convection as detailed in Section 4 have shown that
the resulting neural network can lead to unphysical results, specifically the formation of ice in the lower levels of the atmosphere
and violation of the second law of thermodynamics (negative growth in total entropy for an adiabatic system).

\section{Description of the dynamical model}

The discrete numerical model used so solve \eqref{euler} uses a third order discontinuous Galerkin (DG) spatial 
discretisation, with Gauss-Lobatto-Legendre quadrature and a tensor product Lagrange polynomial basis 
with roots collocated with the quadrature points \cite{Ricardo26} within the discrete discontinuous space 
$S\subset L^2(\Omega)$. The material transport terms are formulated 
so as to provably dissipate tracer variance in the absence of source terms for the piecewise discontinuous 
basis functions \cite{Ricardo24}, with the tracer variance defined as
\begin{equation}
\mathcal{S}_{\gamma} = \frac{1}{2}\int\rho\gamma^2\mathrm{d}\Omega
\end{equation}
for $\gamma\in\{\eta,q_v,q_l,q_i\}$, and the transport terms are cast in material form \cite{Lee26}. In the absence
of non-equilibrium sources terms, tracer variance is provably conserved via the left multiplication of \eqref{euler}
by $\delta\mathcal{S}_{\gamma}/\delta y$, where for example $\delta\mathcal{S}_{\eta}/\delta y = [\boldsymbol{0},\frac{1}{2}\eta^2,\rho\eta,0,0,0]$. 
Conservation
or dissipation of tracer variance in this form has been shown to enhance model stability, since tracer variance
is a convex invariant of the continuous equations in the absence of source terms. 

These tracer variance conserving advective terms are derived by first taking the flux form advection equation for
$\gamma$, $\frac{\partial\rho\gamma}{\partial t} = -\nabla\cdot(\rho\gamma\boldsymbol{u}) = -\frac{1}{2}(\rho\boldsymbol{u}\cdot\nabla\gamma + \gamma\nabla\cdot(\rho\boldsymbol{u}) + \nabla\cdot(\rho\gamma\boldsymbol{u}))$, and then subtracting the tracer weighted continuity equation,
$\gamma\frac{\partial\rho}{\partial t} = -\gamma\nabla\cdot(\rho\boldsymbol{u})$, which for $\rho > 0$ gives
\begin{equation}
\frac{\partial\gamma}{\partial t} = -\frac{1}{2\rho}\Big(\rho\boldsymbol{u}\cdot\nabla\gamma + \nabla\cdot(\rho\boldsymbol{u}\gamma) - 
\gamma\nabla\cdot(\rho\boldsymbol{u})\Big).
\end{equation} 
We then replace the continuous representations of $\rho,\gamma,\boldsymbol{u}$ by their discrete 
analogues $\rho_h,\gamma_h\in S$, $\boldsymbol{u}_h\in\boldsymbol{S}$. Integration with respect
to the test function $\phi_h\in S$ over the domain $\Omega$ (and applying integration by parts to 
the second term on the right hard side) then gives
\begin{multline}\label{adv_disc}
\int\phi_h\frac{\partial\gamma_h}{\partial t}\mathrm{d}\Omega = -\frac{1}{2}\int\phi_h\boldsymbol{u}_h\cdot\nabla\gamma_h -
\gamma_h\rho_h\boldsymbol{u}_h\cdot\nabla\Big(\frac{\phi_h}{\rho_h}\Big) - 
\frac{\phi_h\gamma_h}{\rho_h}\nabla\cdot(\rho_h\boldsymbol{u}_h)\mathrm{d}\Omega + \\
\frac{1}{2}\int\{\phi_h\boldsymbol{u}_h\cdot\hat{\boldsymbol{n}}\}[\gamma_h] -
\{\gamma_h\rho_h\boldsymbol{u}_h\cdot\hat{\boldsymbol{n}}\}[\frac{\phi_h}{\rho_h}] -
\{\frac{\phi_h\gamma_h}{\rho_h}\}[\rho_h\boldsymbol{u}_h\cdot\hat{\boldsymbol{n}}]\mathrm{d}\Gamma - 
\alpha\int\{|\rho_h\boldsymbol{u}_h\cdot\hat{\boldsymbol{n}}|\}[\frac{\phi_h}{\rho_h}][\gamma_h]\mathrm{d}\Gamma,
\end{multline}
where $\Gamma$ denotes the domain of interior element facets, $\{a\} = (a^+ + a^-)/2$, $[a] = a^+-a^-$ denote mean and jump operators
for the discontinuous variables across the left, $-$ and right, $+$ sides of the facets respectively, and 
$\hat{\boldsymbol{n}}$ is the outward facing unit normal vector at the facets. 
The parameter $\alpha$ accounts for the upwinding of fluxes, such 
that setting $\alpha=1$ results in a fully upwinded solution for which tracer variance is dissipated (in the absence 
of source terms), while $\alpha=0$ denotes a centered flux for which tracer variance is conserved.

The discrete form of the continuity equation is given for $\chi_h\in S$ as
\begin{equation}\label{cont_disc}
\int\chi_h\frac{\partial\rho_h}{\partial t}\mathrm{d}\Omega = -\int\chi_h\nabla\cdot(\rho_h\boldsymbol{u}_h)\mathrm{d}\Omega + 
\int\{\chi_h\}[\rho_h\boldsymbol{u}_h\cdot\hat{\boldsymbol{n}}]\mathrm{d}\Gamma,
\end{equation}
with the discrete form of the pressure gradient term, $-\nabla\Phi$ being given as the adjoint of the right hand side, 
$\int\Phi_h\nabla\cdot\boldsymbol{v}_h\mathrm{d}\Omega -
\int\{\Phi_h\}[\boldsymbol{v}_h\cdot\hat{\boldsymbol{n}}]\mathrm{d}\Gamma$, $\forall\boldsymbol{v}_h\in\boldsymbol{S}$.
Setting $\phi_h=\rho_h\gamma_h=\delta\mathcal{S}/\delta\gamma_h$ in \eqref{adv_disc} and 
$\chi_h=\frac{1}{2}\gamma_h^2=\delta\mathcal{S}/\delta\rho_h$ in \eqref{cont_disc} and adding these expressions
then gives
\begin{equation}
\int\rho_h\gamma_h\frac{\partial\gamma_h}{\partial t}\mathrm{d}\Gamma + \int\frac{1}{2}\gamma_h^2\frac{\partial\rho_h}{\partial t} = 
\frac{1}{2}\int\frac{\partial\rho_h\gamma_h^2}{\partial t}\mathrm{d}\Gamma = \frac{\mathrm{d}\mathcal{S}_{\gamma}}{\mathrm{d}t} =
-\alpha\int\{|\rho_h\boldsymbol{u}_h\cdot\hat{\boldsymbol{n}}|\}[\gamma_h]^2\mathrm{d}\Gamma \le 0,
\end{equation}
such that tracer variance is provably bounded semi-discretely in the absence of thermodynamic source terms.

The discrete pressure gradient terms are then formulated as adjoints of the material transport terms so as 
preserve the skew symmetry of the dynamics and conserve energy in the spatial discretisation \cite{Ricardo26}.
For an intrinsic thermodynamic variable, $\psi = \delta u/\delta\gamma$, $\psi\in\{T,\mu_v,\mu_l,\mu_i\}$, 
the corresponding discrete pressure gradient term is given for $\boldsymbol{v}_h\in\boldsymbol{S}$ as
\begin{multline}\label{pgf_disc}
\int\boldsymbol{v}_h\cdot\frac{\partial\boldsymbol{u}_h}{\partial t}\mathrm{d}\Omega = \ldots +
\frac{1}{2}\int{\psi_h\boldsymbol{v}_h}\cdot\nabla\gamma_h -
\gamma_h\boldsymbol{v}_h\cdot\nabla\psi_h - 
\psi_h\gamma_h\nabla\cdot\boldsymbol{v}_h\mathrm{d}\Omega - \\
\frac{1}{2}\int\{\psi_h\boldsymbol{v}_h\cdot\hat{\boldsymbol{n}}\}[\gamma_h] -
\{\gamma_h\boldsymbol{v}_h\cdot\hat{\boldsymbol{n}}\}[\psi_h] -
\{\psi_h\gamma_h\}[\boldsymbol{v}_h\cdot\hat{\boldsymbol{n}}]\mathrm{d}\Gamma +
\alpha\int\{|\boldsymbol{v}_h\cdot\hat{\boldsymbol{n}}|\}[\psi_h][\gamma_h]\mathrm{d}\Gamma.
\end{multline}
Setting $\phi_h = \rho_h\psi_h = \delta\mathcal{H}/\delta\gamma_h$ in \eqref{adv_disc} and 
$\boldsymbol{v}_h=\rho_h\boldsymbol{u}_h = \delta\mathcal{H}/\delta\boldsymbol{u}_h$ in 
\eqref{pgf_disc} leads to the leads to the cancellation of these terms, showing that energy exchanges
are exactly balanced between the advective and pressure gradient forcing terms.

The model uses a stiffly stable third order Runge-Kutta time integrator \cite{SO88}, which unlike the spatial
discretisation does not conserve energy or tracer variance exactly, such that the conservation errors 
decrease with the truncation error of the third order time step.
The phase exchange terms are limited such that these are only applied if the
sum of the advective and phase exchange terms do not result in negative mass fractions. Since the vapour to
liquid $B$, vapour to ice $C$ and liquid to ice $D$ terms are applied separately, they may be limited
separately also. The entropy increment \eqref{entropy_inc} derived from the phase exchanges is consistent 
with the limited form of the phase exchanges.

\section{Results}

In order to assess the non-equilibrium phase exchanges and entropy increment derived from the 
neural network, we apply the DG model to a standard test case for the ascension of a moist non-hydrostatic bubble \cite{Bryan02}, 
extended here to a three phase configuration including ice.
Both the neutral network driven DG model and its initial conditions are entirely independent of the numerical model 
and training data used to learn the neural network. 
The results are then compared to those for a three-phase bubble under the assumption of full thermodynamic equilibrium 
\cite{Ricardo26}, as well as to those from the LFRic model using the CASIM non-equilibrium microphysics scheme. 
For the equilibrium model, the DG solver is the same as for the non-equilibrium 
case, except that the total water is advected as a single prognostic variable, and at each Runge-Kutta sub step the 
partition of the moisture into vapour, liquid and ice is determined via the solution of a minimisation problem where 
the Gibbs free energies of the different phases are matched. This is in contrast to the non-equilibrium DG solver, where 
the three phases are advected independently, and the exchanges derived from the neural network described above.
The LFRic model differs from the neural network driven DG model in that it uses a mixed, lowest order finite element discretisation and
implicit time stepper for the fast wave dynamics, and a flux form semi-Lagrangian finite volume explicit transport scheme, with the
CASIM microphysics applied in order to determine the phase exchanges only once per time step. Moreover the second moment of the 
CASIM microphysics scheme will not have sufficient time
to spin up over the course of the simulation, so this should be regarded as a single moment scheme only, and not a double moment
scheme as is the case for the model simulation used to produce the training data.

The domain is discretised for the DG model as $80\times 80$ $3^{\mathrm{rd}}$ order elements in both dimensions for
$\Omega=10^{4}m\times10^{4}m$ with the time step determined from the acoustic wave speed CFL of 0.5 ($\Delta t=0.05s$) 
and wall boundary conditions of
$\boldsymbol{u}\cdot\hat{\boldsymbol{n}}|_{\partial\Omega} = 0$, where $\partial\Omega$ is the domain boundary.
The fluxes are upwinded by setting $\alpha=1$ in \eqref{adv_disc} such that the variance of all the tracers are dissipated
in the absence of source terms. For the LFRic configuration, the domain is discretised with an equivalent $240\times 240$
lowest order elements in each dimension, but with a larger time step of $\Delta t=2.0s$, owing to the use of the implicit
Crank-Nicolson time stepping scheme, and periodic rather than wall boundary conditions on the horizontal boundaries.

The model is initialised with a 
constant background dry potential temperature of $\theta=300 K$, and with the Exner pressure derived from hydrostatic
balance as $\Pi = 1 - gz/(c_{pd}\theta)$, where $c_{pd}=1004.0J/kg/K$ is the specific heat of dry air at constant pressure.
The mean pressure and density are then initialised respectively as $p = p_{0d}\Pi^{c_{pd}/R_d}$, $\rho = p/(R_d\Pi\theta)$,
where $p_{0d} = 10.0^{5}Pa$ is the reference pressure for dry air at the surface and $R_d=287.0 J/kg/K$ is the specific gas constant 
of dry air. From these the initial vapour profile is determined for a relative humidity of 0.95 (with the initial liquid and ice 
phases set to 0 everywhere, and all mean profiles being constant in the horizontal). A density perturbation 
is then added of the form $\rho' = -\frac{2.0\rho}{\theta}\cos^2(\frac{\pi r}{2r_{max}})$ where $r<r_{max}$ for $r_{max}=2000m$ and
$r = \sqrt{x^2 + (z-r_{max})^2}$.

The entropy perturbation (from its background value of $300K$) and vapour, liquid and ice mass fractions are given for the
non-equilibrium simulation at times $400s, 500s, 600s$ in Figures 2, 3, and 4 respectively. While the entropy is
in units of $J/kg/K$ the mass fractions are all in dimensionless units of $kg/kg$. These show that as the bubble
ascends there is a steady conversion of vapour to liquid in the ascending core of the bubble, while the upper and outer
extents of the bubble are converted into ice. 

\begin{figure}[!hbtp]
\centering
\includegraphics[width=0.80\textwidth,height=0.60\textwidth]{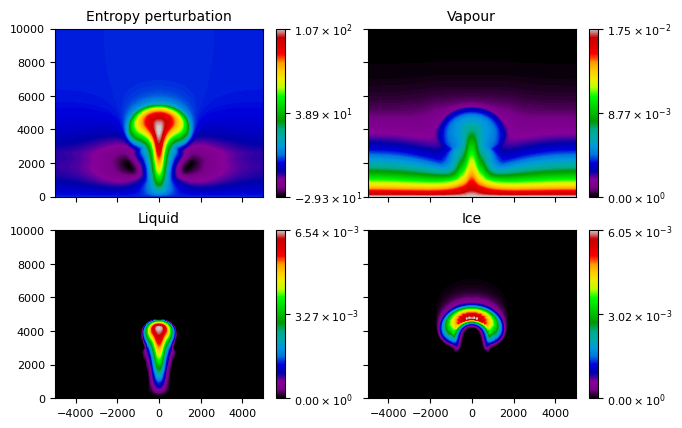}
\caption{Moist bubble with non-equilibrium phase changes, 400.0s. 
}
\end{figure}

\begin{figure}[!hbtp]
\centering
\includegraphics[width=0.80\textwidth,height=0.60\textwidth]{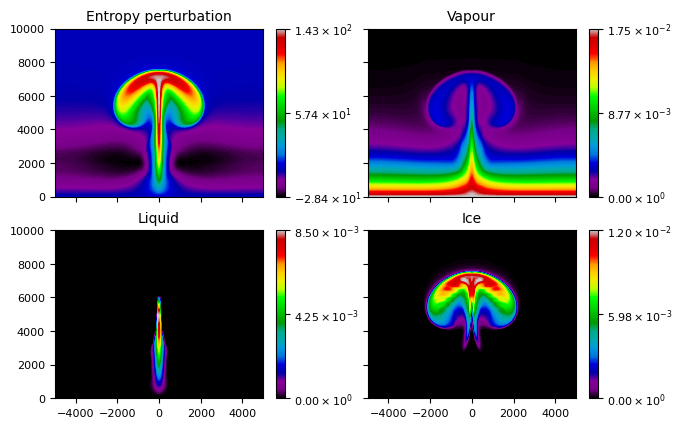}
\caption{Moist bubble with non-equilibrium phase changes, 500.0s}
\end{figure}

\begin{figure}[!hbtp]
\centering
\includegraphics[width=0.80\textwidth,height=0.60\textwidth]{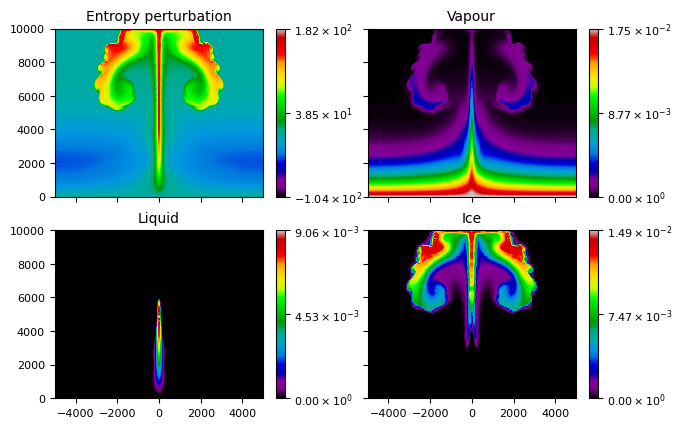}
\caption{Moist bubble with non-equilibrium phase changes, 600.0s}
\end{figure}

The corresponding solutions for the entropy perturbation, vapour, liquid and ice mass fractions for the equilibrium configuration
at times $400s$, $500s$, $600s$ are given in Figures 5, 6, and 7 respectively. 
The results are qualitatively similar to those for 
the non-equilibrium configuration, with ice in the outer extents of the ascending bubble, a liquid core and vapour in the lower levels
of the atmosphere. However the different water phases are more spatially distinct with slightly larger magnitudes and sharper gradients
between phases, perhaps due to some of the logical conditions that are imposed on the equilibrium partition of moisture species,
such as that liquid water is not permitted below the freezing temperature.
The close similarity between the non-equilibrium and equilibrium results is perhaps due to the small time step (0.5 the CFL of the
acoustic wave speed) such that the deviation of the non-equilibrium exchanges from the equilibrium solution is small.

\begin{figure}[!hbtp]
\centering
\includegraphics[width=0.80\textwidth,height=0.60\textwidth]{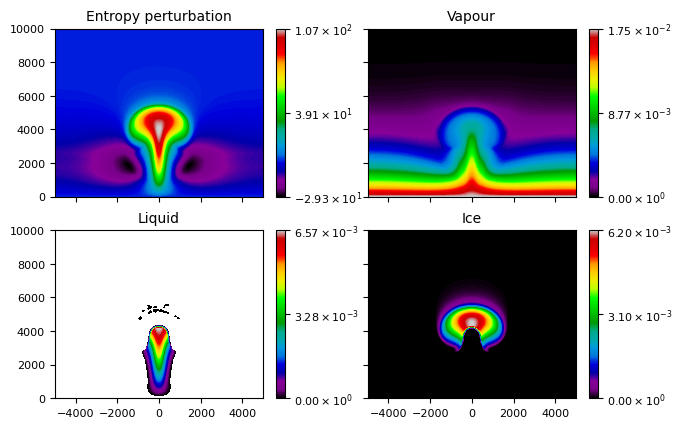}
\caption{Moist bubble with equilibrium phase changes, 400.0s}
\end{figure}

\begin{figure}[!hbtp]
\centering
\includegraphics[width=0.80\textwidth,height=0.60\textwidth]{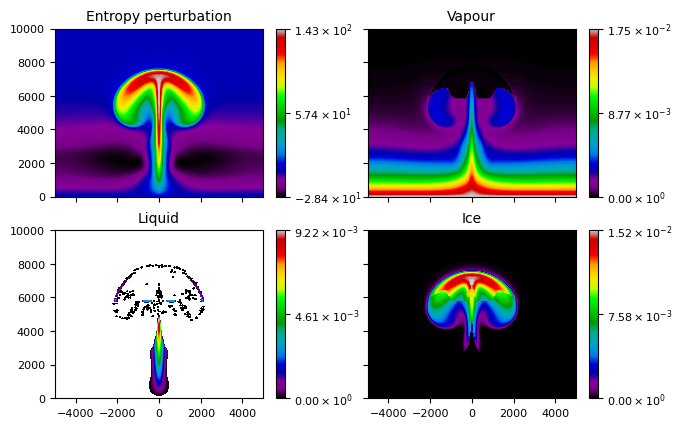}
\caption{Moist bubble with equilibrium phase changes, 500.0s}
\end{figure}

\begin{figure}[!hbtp]
\centering
\includegraphics[width=0.80\textwidth,height=0.60\textwidth]{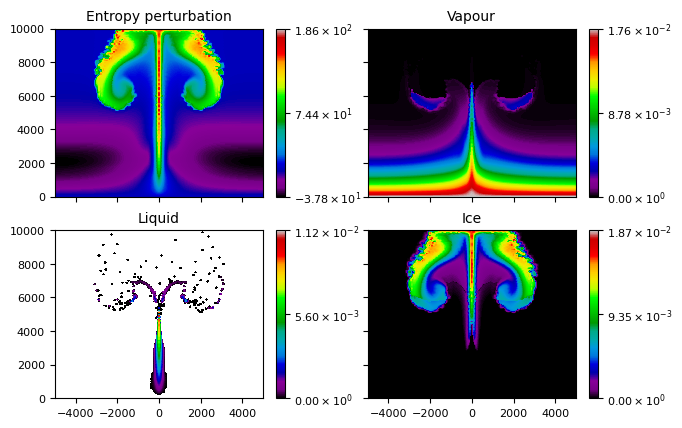}
\caption{Moist bubble with equilibrium phase changes, 600.0s}
\end{figure}

The results for the corresponding LFRic model simulation at a time of $320s$ are presented in Figure 8. Note that while the mixing ratios 
of the different water species as used in LFRic have been converted into mass fractions, the thermal variable presented here is the 
potential temperature perturbation, and not the entropy perturbation as used in the DG model. Also note that the bubble ascends somewhat
faster in the LFRic simulation, with a slightly larger conversion of vapor to ice, with the DG model taking approximately $500s$ to reach
the same height with a more compact shape. Also note that the LFRic simulation sees the conversion of some vapour to liquid in the
lower levels of the atmosphere and to ice above, whereas in the DG simulations the conversions are all in the ascending bubble. The 
extent to which these discrepancies may be attributed to the different representations of the microphysics or the different formulations
of the dynamics is difficult to ascertain. Nevertheless the rates of bubble ascension and magnitudes of ice and liquid formation are
broadly similar between the different implementations. It is also notable that the LFRic configuration fails shortly after this time,
with the bubble approaching the top of the domain, whereas the DG model, which conserves energy and bounds tracer variance, is stable
long after the bubble reaches the top of the domain.

\begin{figure}[!hbtp]
\centering
\includegraphics[width=0.80\textwidth,height=0.60\textwidth]{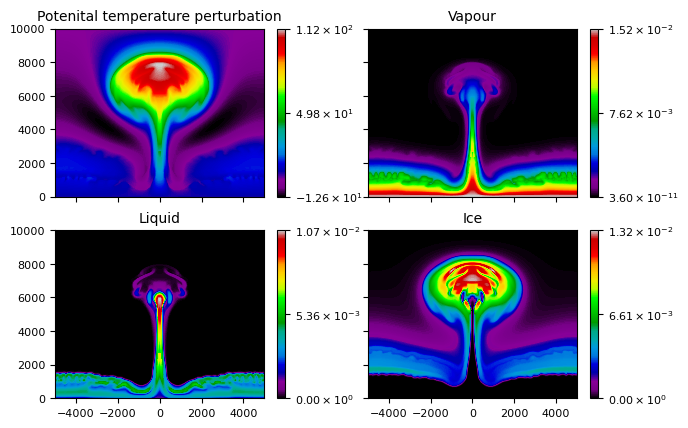}
\caption{Non-equilibrium LFRic configuration of the moist bubble test case, 320.0s}
\end{figure}

This steady conversion from vapour to liquid and ice for the DG model is also observed in the time series of the 
global integrals of the mass of the different phases, $\int\rho q_k\mathrm{d}\Omega$, $k\in\{v,l,i\}$ presented in Figure 9
for both the non-equilibrium and equilibrium configurations. While the magnitudes of these exchanges are comparable for the
two different configurations, the change in the global integral of the entropy, $\int\rho\eta\mathrm{d}\Omega$, is qualitatively different. For the 
non-equilibrium configuration the total entropy grows monotonically, which is consistent with the second law of thermodynamics
for an adiabatic system, while the equilibrium formulation appears to violate this principle with negative entropy growth
for much of the simulation, until the time where the bubble starts to reach the top of the domain.

\begin{figure}[!hbtp]
\centering
\includegraphics[width=0.48\textwidth,height=0.36\textwidth]{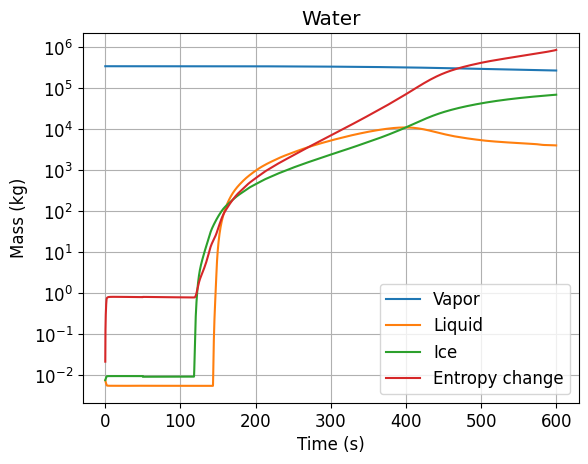}
\includegraphics[width=0.48\textwidth,height=0.36\textwidth]{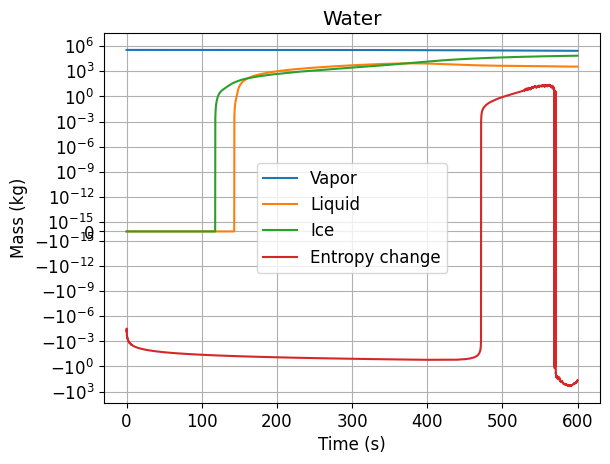}
\caption{Time series of the mass fractions of the different water phases for the non-equilibrium (left)
and equilibrium (right) configurations of the moist bubble.}
\end{figure}

Figure 10 gives the power associated with the change in energy per unit time due to phase exchanges for the non-equilibrium 
configuration. 
These are determined via the left multiplication of the thermodynamic source terms in \eqref{euler} by the vector 
$\delta\mathcal{H}/\delta y$ as
\begin{subequations}
\begin{align}
P_{vapour} &:= \rho^2\mu_v(B(\mu_v-\mu_l) + C(\mu_v - \mu_i)),\\
P_{liquid} &:= \rho^2\mu_l(B(\mu_l-\mu_v) + D(\mu_l - \mu_i)),\\
P_{ice} &:= \rho^2\mu_i(C(\mu_i-\mu_v) + D(\mu_i - \mu_l)),\\
P_{entropy} &:= \rho^2T^2A = -\rho^2((\mu_v-\mu_l)^2{B} + (\mu_v-\mu_i)^2{C} + (\mu_l-\mu_i)^2{D}).
\end{align}
\end{subequations}
The sum of all the power terms (water phases and entropy) is 
$\mathcal{O}(10^{-8}) W$, which is $\mathcal{O}(10^{-15})$ times smaller than the magnitude of the power exchanges, suggesting
that these exchanges are balanced to machine precision, and do not contribute to any change in energy, as determined by 
\eqref{entropy_inc}. 

\begin{figure}[!hbtp]
\centering
\includegraphics[width=0.48\textwidth,height=0.36\textwidth]{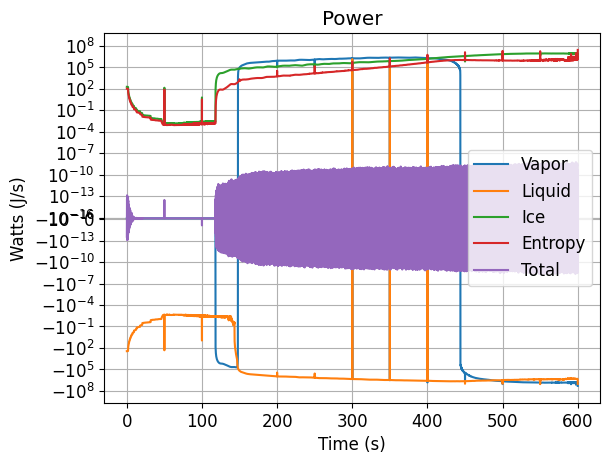}
\caption{Time series of the power associated with changes of energy between phases
for the non-equilibrium moist bubble.}
\end{figure}

\section{Conclusions}

A neural network is trained to learn the non-equilibrium exchanges between vapour, liquid and ice from output of the LFRic 
model using the CASIM microphysics parameterisation scheme for a km-scale regional simulation. This network is 
then implemented within a non-equilibrium discontinuous Galerkin model of the moist compressible Euler equations 
at sub-km resolution, where the separate phases of 
moisture are advected independently and the entropy increment associated with the phase exchanges is derived via conservation
of energy. Results for the simulation of a moist bubble are compared to those for a physics based representation of 
phase exchanges under the assumption of full thermodynamic equilibrium, and to those from the LFRic model using the original
CASIM microphysics scheme from which the neural network training data is generated (albeit with a different model configuration). 

The non-equilibrium and equilibrium DG model configurations are in close agreement, however the non-equilibrium formulation
correctly respects the second law of thermodynamics with positive total entropy growth throughout the simulation, whereas
the equilibrium configuration does not. Both DG configurations show a slightly slower bubble ascent, and slightly lower 
conversion from vapour to liquid and ice than the corresponding LFRic configuration. However it is unclear to what extent
these discrepancies are attributable to the representation of moist processes or to differences in the representation of
dynamics (discretisation, time stepping schemes, boundary conditions).

In future work the neural network could potentially be applied as a substitute for the original microphysics parameterisation
within a convection resolving regional NWP model. In this context the inference of the network may potentially be
more efficient than the existing parameterisation, since the second moment (the particle number) is not required as an 
input to the neural network. It may also potentially be more stable, since energy conservation for moist processes at non-equilibrium
may be assured via the derivation of the entropy increment associated with the phase exchanges. 

In addition to the energy conserving, thermodynamically consistent formulation of non-equilibrium phase exchanges presented here,
it would also be desirable to ensure that thermodynamic exchanges always respect the second law of thermodynamics in the absence 
of diabatic source terms. While this is observed for the results of the neural network generated phase exchanges presented here,
it is not strictly enforced by either the training of the network or the forcings generated from the network in the dynamical solver. 
This is challenging to achieve however, since the
training data derived from the model output implicitly accounts for thermal source terms as well as lateral and bottom 
boundary forcings.
Potentially some form of global penalty term for violations of the second law of thermodynamics, together with a process for 
filtering out the thermodynamic source terms from the training data, may be used to account for this as part of future work.
One might also train the network from data derived from a more idealised model configuration, such as a convection resolving
aquaplanet without orography or boundary layer parameterisations, so as to eliminate any aliasing from these processes into 
the training data.

Finally, we note that the phase exchanges derived from the neural network exhibit a strong sensitivity to the configuration
of the training process, such as the batch size, the weighting of the source terms with respect to the transport terms in the 
loss function and the magnitude of the training loss. For example, 
further reduction of the loss function for the neural network training process presented here results in
unrealistic representations of thermodynamic processes, such as negative entropy 
growth and freezing of water vapour at the surface. 
The neural network was therefore selected based on both the convergence of the loss and physical diagnostics in the independent 
moist bubble test case. 
More extensive explorations into the generation of training data and the application of the network to different physical regimes should 
also be undertaken as part of future work.

\section*{Acknowledgements}

David Lee would like to than Dr. Ian White for a careful reading and insightful feedback of an early draft of this manuscript.

\section*{Appendix A: Conversion of training data}

The LFRic model represents the moisture species as mixing ratios for vapour, $m_v$, cloud liquid, $m_{cl}$, cloud ice, $m_{ci}$,
rain $m_r$, snow, $m_s$ and graupel, $m_g$. These are converted to mass fractions for vapour, liquid and ice as:
\begin{subequations}
\begin{align}
q_v &= \frac{m_v}{1 + m_v + m_{cl} + m_{ci} + m_r + m_s + m_g}, \\
q_l &= \frac{m_{cl} + m_r}{1 + m_v + m_{cl} + m_{ci} + m_r + m_s + m_g}, \\
q_i &= \frac{m_{ci} + m_s + m_g}{1 + m_v + m_{cl} + m_{ci} + m_r + m_s + m_g}.
\end{align}
\end{subequations}
The density is then computed from the potential temperature, $\theta$ and Exner pressure, $\Pi$ as
\begin{equation}
\rho = \frac{p_{0d}}{R_d\theta (1 + m_vR_v/R_d)}\Pi^{c_{vd}/R_d}
\end{equation}
In order to compute the entropy and temperature we first compute the specific heat at constant pressure, $c_p$ and volume, $c_v$ 
and the gas constant for all phases $R$ as
\begin{subequations}
\begin{align}
c_p &= c_{pd}q_d + c_{pv}q_v + c_lq_l + c_iq_i, \\
c_v &= c_{vd}q_d + c_{vv}q_v + c_lq_l + c_iq_i, \\
R &= R_dq_d + R_vq_v,
\end{align}
\end{subequations}
where $c_{pd}, c_{vd}$ are the specific heats of dry air at constant pressure and volume respectively, $c_{pv}, c_{vv}$ are the
specific heats of vapour at constant pressure and volume, $c_l, c_i$ are the specific heats of liquid water and ice respectively 
(recalling that liquid water and ice are assumed to be incompressible), and $R_d, R_v$ are the specific gas constants of dry air and vapour.
For a full list of thermodynamic constants used here see \cite{Bowen22a}. The entropy is then given as
\begin{multline}\label{eq::eta}
\eta = c_p\log\Bigg(\frac{\theta}{T_0}\Bigg) - R\log\Bigg(\frac{p_0}{R}\Bigg) - R_dq_d\log\Bigg(\frac{R_dq_d}{p_{0d}}\Bigg) - R_vq_v\log\Bigg(\frac{R_vq_v}{p_{0v}}\Bigg) + \\
q_d(c_{pd}\log(T_0) - R_d\log(p_{0d})) + q_v\Bigg(c_{pv} + \frac{L_{00s}}{T_0}\Bigg) + q_l\Bigg(c_l + \frac{L_{00f}}{T_0}\Bigg) + q_ic_i,
\end{multline}
where $p_{0d}, p_{0v}$ are the reference pressures, of dry air and vapour, $p_0=p_{0d}+p_{0v}$, $T_0$ is the reference temperature and
$L_{00s}, L_{00f}$ are the specific latent heats of fusion and sublimation.
The last four terms in \eqref{eq::eta} correspond to the reference entropy, $\eta_0$.

Having determined the prognostic variables, $\eta,q_j$, the intrinsic thermodynamic variables may then be determined from the
internal energy 
\begin{equation}
u = T_0c_v\exp\Bigg(\frac{\eta - \eta_0}{c_v}\Bigg)
\Bigg(\frac{q_d\rho}{\rho_{0d}}\Bigg)^{q_dR_d/c_v}
\Bigg(\frac{q_v\rho}{\rho_{0v}}\Bigg)^{q_vR_v/c_v} + q_vL_{00s} + q_lL_{00f},
\end{equation}
where $\rho_{0d}, \rho_{0v}$ are the reference densities of dry air and vapour respectively as
\begin{subequations}
\begin{align}
T &= T_0\exp\Bigg(\frac{\eta - \eta_0}{c_v}\Bigg)
\Bigg(\frac{q_d\rho}{\rho_{0d}}\Bigg)^{q_dR_d/c_v}
\Bigg(\frac{q_v\rho}{\rho_{0v}}\Bigg)^{q_vR_v/c_v},\\
p &= \rho RT,\\
\mu_d &= T\Bigg(R_d\log\Bigg(\frac{q_d\rho}{\rho_{0d}}\Bigg) - c_{vd}\log\Bigg(\frac{T}{T_0}\Bigg) + 
c_{pd} + c_{pd}\log\Bigg(\frac{p_{0d}}{T_0}\Bigg) - c_{vd}\log(p_{0d})\Bigg), \\
\mu_v &= T\Bigg(R_v\log\Bigg(\frac{q_v\rho}{\rho_{0v}}\Bigg) - c_{vv}\log\Bigg(\frac{T}{T_0}\Bigg) - \frac{L_{00s}}{T_0}\Bigg) + L_{00s}, \\
\mu_l &= -T\Bigg(c_{l}\log\Bigg(\frac{T}{T_0}\Bigg) + \frac{L_{00f}}{T_0}\Bigg) + L_{00f}, \\
\mu_i &= -Tc_{i}\log\Bigg(\frac{T}{T_0}\Bigg).
\end{align}
\end{subequations}

\end{document}